\documentclass[11pt]{article}
\usepackage{amsmath, amsthm}
\usepackage{latexsym}
\usepackage{graphicx}
\makeatletter
\newfont{\footsc}{cmcsc10 at 8truept}
\newfont{\footbf}{cmbx10 at 8truept}
\newfont{\footrm}{cmr10 at 10truept}
\makeatother
\newtheorem{theorem}{\bf Theorem}
\newtheorem{proposition}{\bf Proposition}
\newtheorem{corollary}{\bf Corollary}

\newtheorem{conjecture}{\bf Conjecture}
\begin{document}

\title{Bounds on the Location of the Maximum Stirling Numbers of the Second Kind}

\author{Yaming Yu\\
\small Department of Statistics\\[-0.8ex]
\small University of California\\[-0.8ex]
\small Irvine, CA 92697, USA\\[-0.8ex]
\small \texttt{yamingy@uci.edu}}

\date{}
\maketitle
\begin{abstract}
Let $S(n, k)$ denote the Stirling number of the second kind, and let $K_n$ be such that  
$$S(n, K_n-1)<S(n, K_n)\geq S(n, K_n+1).$$
Using a probabilistic argument, we show that, for all $n\geq 2$,
$$\lfloor e^{w(n)}\rfloor -2\leq K_n\leq \lfloor e^{w(n)}\rfloor +1,$$
where $\lfloor x\rfloor$ denotes the integer part of $x$, and $w(n)$ denotes Lambert's W function. 
\end{abstract}

\section{Introduction}
The Stirling number of the second kind, denoted $S(n, k)$, plays a fundamental role in many combinatorial problems.  It counts 
the number of partitions of $\{1, \ldots, n\}$ into $k$ non-empty, pairwise disjoint subsets, and may be defined  recursively 
as 
$$S(n, k)=S(n-1, k-1)+kS(n-1, k),\quad n\geq 1,\ k\geq 1,$$
together with $ S(0, 0)=1,\ S(n, 0)=0,\ n\geq 1.$

According to Harper \cite{H}, for each $n\geq 1$, the polynomial $\sum_{k=0}^n S(n, k)x^k$ has only real zeros.  By Newton's 
inequalities (\cite{HLP}, p. 52), $\log S(n, k)$ is strictly concave in $k$.  It follows that there exists some $1\leq K_n\leq 
n$ such that  
$$S(n, 1)<\ldots <S(n, K_n)\geq S(n, K_n+1)>\ldots>S(n, n).$$
In other words, the sequence $S(n, k),\ k=1,\ldots, n,$ is unimodal, $K_n$ being a unique mode if $S(n, K_n)\neq S(n, 
K_n+1)$.

Determining the value of $K_n$ is an old problem (\cite{K1, K2, Harb, B, Dob, K3, W, M, C}).  A related long-standing 
conjecture (\cite{W, CP, KMR}) is that there exists no $n>2$ such that $S(n, K_n)=S(n, K_n+1)$.  See \cite{CP} for a 
historical sketch and recent developments.

In particular, Canfield and Pomerance \cite{CP} noted that 
\begin{equation}
\label{kn}
K_n\in \{\lfloor e^{w(n)}\rfloor -1,\ \lfloor e^{w(n)}\rfloor \}
\end{equation}
for both $2\leq n\leq 1200$ and $n$ large enough (no specific bound is known on how large $n$ has to be; see also \cite{C}).  
Here and in what follows, $\lfloor x\rfloor$ denotes the integer part of $x$ and $w(n)$ is Lambert's W function defined by 
$$n=w(n)e^{w(n)}.$$
Based on this, it seems likely that (\ref{kn}) holds for all $n$.  The purpose of this note is to present the following 
non-asymptotic bounds.
\begin{theorem}
\label{m}
\begin{equation}
\label{main}
\lfloor e^{w(n)}\rfloor -2\leq K_n\leq \lfloor e^{w(n)}\rfloor +1,\quad n\geq 2.
\end{equation}
\end{theorem}

Theorem \ref{m} can be compared with the non-asymptotic bounds of Wegner \cite{W}:
\begin{align}
\label{upper}
K_n &< \frac{n}{\log n-\log\log n}, \quad n\geq 3;\\
\label{lower}
K_n &> \frac{n}{\log n}\left(1+\frac{\log\log n-1}{\log n}\right), \quad n\geq 31.
\end{align}
Note that the upper and lower bounds in (\ref{main}) differ by 3, whereas the difference between the upper bound
(\ref{upper}) and the lower bound (\ref{lower}) tends to $\infty$ as $n\rightarrow\infty$.  More precisely, it can be shown
(details omitted) that the upper bound in (\ref{main}) implies (\ref{upper}) if $n\geq 7$, and the lower bound in
(\ref{main}) implies (\ref{lower}) if $n\geq 34$.

In Section 2 we prove (\ref{main}) using a probablistic result of Darroch \cite{D}.  The possibility of further refinements 
is discussed in Section 3.

\section{Proof of (\ref{main})}
Recall Dobinski's formula
\begin{equation}
\label{dob0}
e^x \sum_{k=1}^n S(n, k) x^k = \sum_{k=1}^\infty \frac{k^n x^k}{k!},\quad n\geq 1.
\end{equation}
In particular
\begin{equation}
\label{dob1}
e\sum_{k=1}^n S(n, k)=\sum_{k=1}^\infty \frac{k^n}{k!}.
\end{equation}
Dividing (\ref{dob0}) by (\ref{dob1}) we get 
$$\left(\sum_{k=0}^\infty \frac{1}{e k!} x^k\right) \sum_{k=1}^n \frac{S(n, k)}{\sum_{i=1}^n S(n, i)} x^k = \sum_{k=1}^\infty 
\frac{k^n /k!}{\sum_{i=1}^\infty i^n/i!} x^k.$$
This has the following interpretation.  If we let $S$ be a random variable with probability mass function (pmf) $\Pr(S=k)=S(n, 
k)/\sum_{i=1}^n S(n, i),\ k=1, \ldots, n,$ and let $Z$ be a ${\rm Poisson}(1)$ random variable independent of $S$, then the 
pmf of $S+Z$ is 
$$\Pr(S+Z=k)=\frac{k^n/k!}{\sum_{i=1}^\infty i^n/i!},\quad k=1, 2, \ldots$$

While the mode of $S$ is hard to determine, that of $S+Z$ is straightforward.  (As usual, we call a random variable $X$ 
on $\{0, 1, \ldots\}$ unimodal if its pmf is unimodal, and call any mode of the pmf a mode of $X$.)  To relate the mode of $S$ 
to that of $S+Z$, we invoke a classical result of Darroch \cite{D} (see Pitman's survey \cite{P}).  Note that $S$ 
can be written as a sum of $n$ independent Bernoulli random variables since the polynomial $\sum_{k=1}^n S(n, k) x^k$ has 
only real zeros.

\begin{theorem}[\cite{D}]
\label{thmd}
Let $X_i,\ i=1, \ldots, n,$ be independent Bernoulli random variables, i.e., each $X_i$ takes values on $\{0, 1\}$.  Then for 
any mode $m$ of $S=\sum_{i=1}^n X_i$
$$|m-ES|<1.$$
\end{theorem}

As a consequence of Theorem \ref{thmd}, we have 
\begin{proposition}
\label{prop1}
Let $S=\sum_{i=1}^n X_i$ be a sum of independent Bernoulli random variables.  Let $Z$ be a ${\rm Poisson}(1)$ random variable 
independent of $S$.  Assume $S+Z$ has a unique mode $m_1$, and denote any mode of $S$ by $m_0$.  Then 
\begin{equation}
\label{m01}
m_0\leq m_1\leq m_0+2.
\end{equation}
\end{proposition}
{\bf Proof.}  Note that, since the pmfs of $S$ and $Z$ are both log-concave, the pmf of $S+Z$ is log-concave
and hence unimodal.  Denote $\mu=ES$.  By Darroch's rule, $|\mu-m_0|<1$.  We show that Darroch's rule applies to $S+Z$, i.e., 
$|\mu+1-m_1|<1$.  The claim then readily follows.  Let $Z_k,\ k\geq 2,$ be ${\rm Binomial}(k, 1/k)$ random variables, 
independent of $S$.  Then $S+Z_k$ is a sum of independent Bernoullis for which Darroch's rule applies; if we let $m_k$ be a 
mode of $S+Z_k$, then $|\mu+1-m_k|<1$.  Moreover, assuming $m_1$ is the unique mode of $S+Z$, we have 
$\lim_{k\rightarrow\infty} m_k=m_1$.  Thus $|\mu+1-m_1|<1$. \qed

On the other hand, we have 
\begin{proposition}
\label{prop2}
For $n\geq 2$, the sequence $k^n/k!,\ k=1, 2,\ldots,$ is unimodal with a unique mode at either $k=\lfloor e^{w(n)}\rfloor$ or 
$k=\lfloor e^{w(n)}\rfloor+1$.
\end{proposition}
{\bf Proof.}  Denote $u= e^{w(n)}$ and consider the ratio
$$f(k)=\frac{(k+1)^n/(k+1)!}{k^n/k!}=\frac{(k+1)^{n-1}}{k^n}.$$
It is easy to see that $f(k)\neq 1$ for all $k\geq 1$.  We also show that $f(k)>1$ if $k<u-1$ (i.e., $k\leq \lfloor 
u\rfloor-1$) and $f(k)<1$ for $k>u$ (i.e., $k\geq \lfloor u\rfloor +1$).  The claim then follows.

Noting that $f(k)$ decreases in $k$, we only need to show $f(u-1)>1$ and $f(u)<1$.  However, direct calculation gives 
\begin{align*}
\log f(u-1) &=-w(n)-n\log\left(1-e^{-w(n)}\right)\\
            &>-w(n)-n\left(-e^{-w(n)}\right)=0;\\
\log f(u)   &=n\log\left(1+e^{-w(n)}\right)-\log\left(e^{w(n)}+1\right)\\
            &<ne^{-w(n)}-\log\left( e^{w(n)}\right)=0. \qed
\end{align*}

Then we obtain (\ref{main}) as a consequence of Propositions \ref{prop1} and \ref{prop2}.
\begin{corollary}
\label{coro}
Let $n\geq 2$, and denote $k_*=\lfloor e^{w(n)}\rfloor$.  If $k_*^n/k_*!>(k_*+1)^n/(k_*+1)!$, then $k_*-2\leq K_n\leq k_*$; 
otherwise $k_*-1\leq K_n\leq k_*+1$.  At any rate (\ref{main}) holds.
\end{corollary} 

\section{Discussion}
A natural question is whether Corollary \ref{coro} can be further improved using this argument.  This leads to an 
investigation of the 
bounds in (\ref{m01}).  It turns out that the lower bound in (\ref{m01}) is achievable.  For example, in the setting of 
Proposition \ref{prop1}, if we let $n=2$ and $\Pr(X_i=1)=1-\Pr(X_i=0)=p_i,\ i=1, 2$, with $p_1=1/3$ and $p_2=2/5$, then  
$m_0=m_1=1$ by direct calculation.  It seems difficult, however, to find an example where the upper bound in (\ref{m01}) is 
achieved.  After some experimentation we suspect that this upper bound is not achievable.  This is further supported by the 
fact that, in the setting of Proposition \ref{prop1}, we always have $m_1\leq m_0+1$ when $n\leq 5$.  To show this, let 
$c_i=\Pr(S=i),\ i=0, 1,\ldots$.  By Newton's inequalities
$$c_{i+1}^2\geq \frac{(i+2)(n-i)}{(i+1)(n-i-1)} c_ic_{i+2},\quad 0\leq i\leq n-2.$$
When $n\leq 5$ and $0\leq i\leq n-2$ we have 
$$\frac{(i+2)(n-i)}{(i+1)(n-i-1)}\geq 2.$$
Thus $c_{i+1}^2\geq 2c_ic_{i+2}$ and $c_{m_0+1}^2\geq 2c_{m_0} c_{m_0+2}$ in particular.  Since $m_0$ is a mode of $S$, we 
know $c_{m_0}\geq c_{m_0+1}$.  Thus
$$c_{m_0}\geq 2c_{m_0+2}.$$
However, a simple calculation gives 
\begin{align*}
e[\Pr(S+Z=m_0+1)-\Pr(S+Z=m_0+2)] &=\sum_{k=0}^{m_0} \frac{c_k}{(m_0-k)!(m_0+2-k)} -c_{m_0+2}\\
&\geq \frac{c_{m_0}}{2}-c_{m_0+2}\geq 0,
\end{align*}
which rules out $m_1= m_0+2$ under the assumption that $m_1$ is the unique mode of $S+Z$. 

\begin{conjecture}
\label{conj}
In the setting of Proposition \ref{prop1}, $m_0\leq m_1\leq m_0+1$.
\end{conjecture}
It is clear that Conjecture \ref{conj} implies a sharper version of (\ref{main})
$$\lfloor e^{w(n)}\rfloor-1\leq K_n\leq \lfloor e^{w(n)}\rfloor+1;$$
this is tantalizingly close to proving (\ref{kn}) for all $n\geq 2$.

\section*{Acknowledgement}
The author would like to thank the anonymous reviewers for their careful reading of the manuscript.

\end{document}